\documentclass[12pt,a4paper,reqno]{amsart}
\usepackage{amsmath}
\usepackage{amsfonts}
\usepackage{amssymb}

\numberwithin{equation}{section}

\newtheorem*{main}{Main Theorem}
\newtheorem*{corollary}{Corollary}%[section]
\newtheorem{lemma}{Lemma}[section]

\def\N{{\mathcal N}}
\def\O{{\mathcal O}}

\def\p{{\mathfrak p}}

\def\q{{\mathfrak q}}
\def\R{{\mathbb R}}

\def\Z{{\mathbb Z}}

\DeclareMathOperator{\Br}{Br}

\DeclareMathOperator{\Gal}{Gal}

\DeclareMathOperator{\sgn}{sgn}

\DeclareMathOperator{\ord}{ord}

\DeclareMathOperator{\Pic}{Pic}

\begin{document}
\title{On the $m$-torsion subgroup of the Brauer group of a global field}
\author[Chi]{Wen-Chen Chi$^*$}
\address{Department of Mathematics\\
National Taiwan Normal University\\
88, Sec.4, Ting-Chou Road\\
Taipei, Taiwan 116\\
E-mail: wchi@math.ntnu.edu.tw}
\author[Liao]{Hung-Min Liao}
\address{Department of Mathematics\\
National Taiwan Normal University\\
88, Sec.4, Ting-Chou Road\\
Taipei, Taiwan 116\\
E-mail: 89340001@ntnu.edu.tw}
\author[Tan]{Ki-Seng Tan$^{\dag}$}
\address{Department of Mathematics\\
National Taiwan University\\
1, Sec.4, Roosevelt Road\\
Taipei, Taiwan 106\\
E-mail: tan@math.ntu.edu.tw}

\thanks{$^*$ The author was supported in part by the National Science Council of Taiwan,
NSC95-2115-M-003-005.
$^{\dag}$ The author was supported in part by the National
Science Council of Taiwan,
NSC95-2115-M-002-017-MY2}

%\subjclass[2000]{11S40 (primary), 11R42, 11R58 (secondary)}

%\keywords{Brumer-Stark Conjecture, Conjecture of Gross, class
%numbers, Iwasawa theory, local Leopoldt conjecture, regulators,
%Stickelberger element, special values of L-functions.}
\maketitle

\begin{abstract}
In this note, we give a short proof of the existence of certain
abelian extension over a given global field $K$. This result implies
that for every positive integer $m$,
there exists an abelian extension $L/K$ of exponent $m$ such that
the $m$-torsion subgroup of $\Br(K)$ equals $\Br(L/K)$.
\end{abstract}

\begin{section}*{}
Assume that a global field K, a
prime number $p$ and a natural number $n$ are given.
We say that a pro-$p$ abelian extension $M/K$ is granted if the
exponent of $M/K$, by which we mean that of $\Gal(M/K)$, divides $p^n$ and for every finite place
$w$ of $M$ the degree $[M_w:K_v]$, $w|v$, also divides $p^n$. A finite place $v$
is said to be supported under $M/K$ if for each $w|v$,
$[M_w:K_v]=p^n$; an archimedean place $v$ is supported under
$M/K$ if $M_w$ is complex for every $w|v$.
The purpose of this note is to give a simple proof of the following.
\begin{main}
Suppose $K$, $p$ and $n$ are given as above. Then there exists a granted extension
that supports every place of $K$.
\end{main}
We should remind the reader that \cite{ks1,ks2,p} put together already gave a
complete proof of the theorem. However, our proof is much simpler.
Furthermore, as these papers have pointed out,
using the structure
theorem of Brauer groups over global fields,
one can deduce from the main theorem the following corollary which gives an affirmative
answer to the question
raised in \cite{as}. Recall that the relative Brauer group $\Br(L/K)$ is the kernel of the
restriction map $\Br(K)\longrightarrow \Br(L)$ (\cite{as}).
\begin{corollary}
For a positive integer $m$,
there always exists an abelian extension $L/K$ of exponent $m$ such that
the $m$-torsion subgroup of $\Br(K)$ equals the relative Brauer group
$\Br(L/K)$.
\end{corollary}
We shall prove the theorem by constructing a sequence
$$M_0\subset M_1
\subset \cdots \subset M_k \subset \cdots$$
of granted extensions over
$K$ and show that $L=\bigcup_{k=0}^{\infty} M_k$ enjoys the desired
property. We will frequently use the fact that if an abelian extension
$K'/K$ is of exponent $p^n$ and is unramified at $v$ then
the degree $[K'_w:K_v]$, $w\mid v$, divides $p^n$. Thus,
in order to check if $K'/K$ is granted, it is sufficient to check the local degrees
at the ramified places.

Consider $\Pic_K :=K^* \backslash\mathbb{A}_K^*/\prod_v
\mathcal{O}_v^*$, where $\mathbb{A}_K^*$ is the ideles group, $v$ runs
through all places of $K$ and for an Archimedean place $v$ we let
$\mathcal{O}_v^*=K_v^*$. For a finite place $\p$, let
$[\p]$ be the
image of any prime element $\pi_{\p} \in K_{\p}^*$ under the natural
map $K_{\p}^* \rightarrow \mathbb{A}_K^* \rightarrow \Pic_K$. The abelian group
$\Pic_K$ is finitely generated and by the Chabotarev density theorem we can
find $\p_1,\p_2,\ldots,\p_\ell$ outside any given finite set of
places such that $[\p_1], [\p_2],\ldots , [\p_\ell]$ form a set of generators of
$\Pic_K$. Denote $S=\{ \p_1,\p_2, \ldots , \p_\ell \}$. For later usage, we choose for every idele
$x \in \mathbb{A}_K^*$ a global element $f_x\in K^*$, which
is unique up to
$\mathcal{O}_S^*$, such that
\begin{equation}\label{e:xyfu}
x=y\cdot  f_x \cdot u
\end{equation}
for some $y\in
\prod_{i=1}^{\ell} K_{\p_i}^*$ and $u \in \prod_v
\mathcal{O}_v^*$.

Put $\Gamma_{-1}=\Pic_K/(p^n\Pic_K)=K^*
\backslash\mathbb{A}_K^*/(\prod_v \mathcal{O}_v^* \cdot(\mathbb{A}_K^*)^{p^n})$.
Class Field Theory (\cite{at})
identifies $\Gamma_{-1}$ as the Galois group of an abelian
extension, denoted as $M_{-1}/K$, which is everywhere unramified.
This together with the obvious fact that the exponent of $\Gamma_{-1}$ divides $p^n$
implies that the extension $M_{-1}/K$ is granted.

Suppose $T$ is a finite set of places with $T\cap
S=\emptyset$ and $\N \subset \prod_{v\in T}\O_v^*$ is an open
subgroup which contains $\O_S^*$ via the natural embedding
$K^* \hookrightarrow\prod_{v\in T}K_v^*$. Set
$$\Gamma(\N):=K^*
\backslash \mathcal{A}_K^* / ((\prod_{v\notin T}\O_v^* \times \N) \cdot
(\mathbb{A}_K^*)^{p^n}).$$
This group is also of exponent dividing $p^n$.
Again, Class Field Theory identifies $\Gamma(\N)$ as the Galois group of an
abelian extension which we denote as $M(\N)/K$.

\begin{lemma}\label{l:1} Let $\N$ be as above.
Then we have the exact sequence
\begin{equation}\label{e:e}
0 \rightarrow \prod_{v\in T}\O_v^*/(\N \cdot \prod_{v\in T}(\O_v^*)^{p^n})
\stackrel{i}{\rightarrow}\Gamma(\N) \stackrel{q}{\rightarrow} \Gamma_{-1}
\rightarrow 0,
\end{equation}
where $i$ is induced from the natural map $\prod_{v\in T}K_v^*
\rightarrow \mathbb{A}_K^*$ and $q$ is the natural quotient map.
\end{lemma}
\begin{proof}
It is enough to show the injectivity of $i$.
Suppose ${\bar z}\in\ker(i)$ is obtained from an element
$z \in \prod_{v\in T} \O_v^*$. Then
$z=f\cdot t\cdot x^{p^n}$ for some $f\in K^*$, $t\in \prod_{v\notin T}\O_v^*\times\N$
and $x\in \mathbb{A}_K^*$. Write $x=y\cdot f_x\cdot u$ as in (\ref{e:xyfu}).
Then we see that $f_x^{p^n}\cdot f \in \O_S^*$. Therefore if $t_T$ and $u_T$
are the $T$-components of $t$ and $u$, then
$z=(f_x^{p^n}\cdot f)\cdot
t_T\cdot u_T^{p^n}\in
\N\cdot (\prod_{v\in T}\O_v^*)^{p^n}$. Hence $\bar z=0$.
\end{proof}

Suppose $p=2$ and
$\infty_1,\ldots ,\infty_s$, $s>0$, are all the real places of $K$.
Put
$\R_j=K_{\infty_j}$ and let $\sgn_j^\flat$ be the sign map $\R_j^* \rightarrow
\R_j^*/\R_{j,+}\simeq \Z/2\Z$. As $M_{-1}/K$ is unramified everywhere, for
each $j$, we can choose a real place of $M_{-1}$ sitting over $\infty_j$ and
use it to define the sign map
$$\sgn_j:M_{-1}^*\longrightarrow \R_j^*\stackrel{\sgn_j^{\flat}}
{\longrightarrow}\Z/2\Z.$$
Consider $U=\{f\in\O_S^*\subset K^* | \sum_{j=1}^{s}\sgn_j(f)\equiv 0 \pmod{2}\}$.
Then either $\O_S^*=U$ or $\O_S^*=U\coprod gU$ for some $g$.
Let $M'_{-1}$ be the
composite of all $M_{-1}(\sqrt{f})$, $f\in U$.
By the Chabotarev density
theorem we can find a finite place $\p$ splitting completely under
$M'_{-1}/K$.
Since on $(M_{-1}^*)^2$ all
the values of the map $\sgn:=\sum_{j=1}^{s}\sgn_j$ equal $0\in\Z/2\Z$, if $g$ exists then
$U\cdot (M_{-1}^*)^2\varsubsetneq \O_S^*\cdot (M_{-1}^*)^2$. In this case, Kummer's Theory
says that
$M'_{-1}(\sqrt{g})/M'_{-1}$ is a
quadratic extension and we choose $\p$ so that it does not
completely split under $M'_{-1}(\sqrt{g})/K$.

Let $\sgn_{\p}^{\flat}$ be the
map $\O_{\p}^* \rightarrow \O_{\p}^*/(\O_{\p}^*)^2 \simeq \Z/2\Z$, for
$(x_{\p},x_1,\ldots,x_s) \in \O_{\p}^* \times \prod_{i=1}^{s}
\R_i^* $ define
$\sgn^{\flat}(x_{\p},x_1,\ldots,x_s)=\sgn_{\p}^{\flat}(x_{\p})+
\sum_{i=1}^s \sgn_i^{\flat}(x_i)$ and
put
$$
\N_0 =\{(x_{\p},x_1,\ldots,x_s) \in \O_{\p}^* \times \prod_{i=1}^{s}
\R_i^* \quad |\quad \sgn^{\flat}(x_{\p},x_1,\ldots,x_s) \equiv
0  \pmod{2} \}.
$$
The extension $M(\N_0)/K$ is unramified outside $\{\p, \infty_1,
\ldots, \infty_s \}$. Since each $\R_i^*\not\subset \N_0$, the
extension is ramified at each $\infty_i$. The map $\sgn^{\flat}$ induces an
isomorphism
$
 \O_{\p}^* \times \prod_{i=1}^s\R_i^* / \N_0
%\times(\O_{\p}^* \times \prod_{i=1}^s\R_i^*)^{2^n}
\simeq \Z/2\Z
$
and from Lemma \ref{l:1} we see that (since $\p$ splits completely under
$M'_{-1}/K$) the local extension of $M(\N_0)/K$ at $\p$ is a
quadratic extension. Therefore $M(\N_0)/K$ is granted and it
supports every real place. We put $M_0=M(\N_0)$ if $p=2$ and $K$
has a real place. Otherwise, put $M_0=M_{-1}$.

We set all finite places of $K$ into a sequence
$\q_1,\q_2,\ldots,\q_k,\ldots$. And we shall construct $M_k$ so that it
supports $\q_k$. Assume that $M_{k-1}$ is already constructed. If
it supports $\q_k$, then we set $M_k=M_{k-1}$. Otherwise, assume
that the decomposition subgroup of
$%\Gamma_{k-1}:=
\Gal(M_{k-1}/K)$ at $\q_k$ is of order $p^m$
with $m<n$. We have $m\geq m_1 + m_2$ where $p^{m_1}$ is the order
of the inertia subgroup of $\Gal(M_{k-1}/K)$ and $p^{m_2}$ is the degree of the residue
field extension of $M_{-1}/K$ at $\q_k$.

Let $S_1=\{v_1,\ldots, v_r\}$ be a finite set of places of $K$ such that
$\q_k\notin S_1$ and $M_{k-1}/K$ is unramified outside
$S_1\cup \{\q_k\}$. We can assume that $S$ is chosen so that
$S\cap (S_1\cup \{\q_k\})=\emptyset$. And
denote
$S_2=S\cup S_1$, $S_3=S_2 \cup \{\q_k\}$. Choose a prime element
$\pi$ at $\q_k$ and let $f_{\pi}\in K^*$ be the global element chosen before.
Then the $S_3$-units group $\O_{S_3}^*$ is the direct product
of $\O_{S_2}^*$ and the infinite cyclic subgroup generated by $f_{\pi}$.

To construct $M_k$, we will need to find a finite place $w$ outside $S_3$,
which splits completely
under $M_{k-1}/K$, and an open subgroup $\N_w\subset\O_w^*$, which contains
$\O_{S_2}^*$ via the natural embedding $K^*\hookrightarrow K_w^*$, such that
the quotient $\O_w^*/\N_w$ is isomorphic to $\Z/p^n\Z$ and contains a subgroup of order $p^{n-m_1}$
generated by
the element $f_{\pi}\pmod{\N_w}$.

%Then we form the extension $M(\N_w)$.

\begin{lemma}\label{l:2}
If $\N_w$ is as above, then the following hold:
\begin{enumerate}
  \item The extension $M(\N_w)/K$ is unramified outside $\{w\}$.
  \item Every finite place $\q\in S_2$ splits completely under
  the extension $M(\N_w)/M_{-1}$.
  \item If $\q\mid\q_k$, then the degree $[M(\N_w)_{\q}:K_{\q_k}]=p^{n-m_1}$.
  \item The extension is totally ramified at $w$ with
  $[M(\N_w)_w:K_w]=p^n$.
\end{enumerate}
\end{lemma}
\begin{proof}
Statements (1) and (4) are from the exact sequence (\ref{e:e}), since
the place $w$ splits completely under $M_{-1}/K$.
To prove (2), let $\pi_{\q}$ be a prime element at $\q$ and
assume that in $\Gamma_{-1}$ the Frobenius at $\q$ is annihilated
by $p^{\mu}$. Then $\pi_{\q}^{p^{\mu}}=f\cdot t\cdot x^{p^n}$ for some $f\in K^*$,
$t\in \prod_v \O_v^*$ and $x\in \mathbb{A}_K^* $. This and the equation
(\ref{e:xyfu}) imply that $f\cdot f_x^{p^n}\in \O_{S_2}^*$. They also imply that
$t_w$, the $w$-component of $t$, is contained in $\N_w$, since
$\O_{S_2}^*\cdot (\O_w^*)^{p^n}\subset \N_w$. Therefore in
$\Gamma(\N_w)$ the Frobenius at $\q$ is also annihilated by
$p^{\mu}$. And (2) is proved.

Let $F_r \in \Gamma(\N_w)$ be the
Frobenius element at $\q_k$. Using the exact sequence (\ref{e:e}), we deduce
that  $p^{m_2}F_r \in \O_w^*/\N_w$ and hence $\pi^{m_2}=f\cdot t\cdot x^{p^n}$
for some $f\in
K^*$, $t\in \prod_v \O_v^*$ and $x\in \mathcal{A}_K^*$. This means
$p^{m_2}F_r$ equals the residue class $t_w\pmod{\N_w}$ and, together with the equation
(\ref{e:xyfu}), this also shows that $t_w\equiv f^{-1}\pmod{(\O_w^*)^{p^n}}$.
We then write
$\pi = z\cdot f_{\pi}\cdot u'$ with $z
\in \prod_{i=1}^sK_{p_i}^*$ and $ u' \in \prod_v\O_v^*$
and use it to deduce that $f^{-1}\cdot f_x^{-p^n}\cdot f_{\pi}^{p^{m_2}} \in \O_{S_2}^*$.
Again, since
$\O_{S_2}^*\cdot (\O_w^*)^{p^n}\subset \N_w$, we have $t_w\equiv f_{\pi}^{p^{m_2}}\pmod{\N_w}$.
Therefore in the quotient $\O_w/\N_w$ the order of the residue class
$t_w\pmod{\N_w}$ equals $p^{n-m_1-m_2}$. And this proves that $F_r$ is of order $p^{n-m_1}$.
\end{proof}

We then set $M_k=M_{k-1}M(\N_w)$. It is clear that the exponent of $\Gal(M_k/K)$ is
$p^n$. The extension $M_k/K$ is unramified outside $S_1\cup\{q_k,w\}$.
Lemma \ref{l:2} implies that $[M_{k,\p}:K_{\q_k}]=p^n$ for $\p\mid \q_k$,
$[M_w:K_w]=p^n$ and if $\q\in S_1$, then $[M_{k,\q}:K_{\q}]=[M_{k-1,\q}:K_{\q}]$.
Therefore $M_k/K$ is granted and it supports $\q_k$.

%By Lemma \ref{l:1} and Lemma \ref{l:2}, it is now clear that we have proved the
%existence of a granted extension $L/K$ which supports every place
%of $K$.

To complete the proof, we need to find $w$ and $\N_w$.
Let us first
consider the case where $\text{char}.(K) \neq p$. Let $\q$ be a place of $M_{k-1}$
sitting over $\q_k$. Then we have $\ord_{\q}(f_{\pi})=p^{m_1}\cdot\ord_{\q_k}(f_{\pi})=p^{m_1}$.
Since $\ord_{\q}(f)=0$ for all $f\in\O_{S_2}^*$, the index $|\O_{S_2}^*\cdot (M_{k-1}^*)^{p^n}:
\O_{S_3}^*\cdot (M_{k-1}^*)^{p^n}|$ is a multiple of $p^{n-m_1}$.
For a global element
$f\in M_{k-1}^*$ let $M_{k-1}(\sqrt[p^n]{f})$ be the Kummer extension
generated by all $p^n$th root of $f$. Let $K_2'$ be the
composite
of all $M_{k-1}(\sqrt[p^n]{f})$, $f \in \O_{S_2}^*$ and let
$K_3'$ be the
composite
of all $M_{k-1}(\sqrt[p^n]{f})$, $f \in \O_{S_3}^*$. Kummer's theory
tells us that the degree $[K_3':K_2']$ is a multiple of $p^{n-m_1}$.
Choose an element $\sigma\in\Gal(K_3'/K_2')\subset \Gal(K_3'/K)$ of order $p^{n-m-1}$
and apply the Chabotarev density theorem to choose
$w$ to be a unramified place under $K_3'/K$ such that the Frobenius at $w$ equals
$\sigma$. We choose $w$ outside $S_3$ and put $\N_w=(\O_w^*)^{p^n}$. It is obvious that
$w$ splits completely under $M_{k-1}$, $\O_{S_2}^*\subset \N_w$ and in the
quotient $\O_w^*/\N_w$ the order of $f_{\pi}\pmod{\N_w}$ equals $p^{n-m_1}$.
Since $K_2'$ contains all the primitive roots of $1$ and $w$ splits completely under
$K_2'/K$, the local ring $\O_w$ also contains all the  primitive roots of $1$.
This implies that $\O_v^*/\N_w \simeq \Z/p^n\Z$. Thus $w$ and $\N_w$ satisfy all
the required conditions.

Finally, we consider the case
where $\text{char}.(K)=p$. Apply the Chabotarev density theorem and choose $w$ to be a
place outside $S_2$, splitting completely under
$M_{k-1}/K$. Recall that the $p$-part $\O_1 = 1+\pi_w \O_w$ of
$\O_w^*$ is the direct product of countable many copies of
$\Z_p$ (\cite{w}) and the local Leopoldt Conjecture holds (\cite{k93})
in the sense
that $\Z_p \otimes_\Z \O_{S_3}^*$ form a direct
summand of $\O_1$. In other words, we have
$$\O_1=\Z_p\otimes \O_{S_2}^*\times \Z_p\otimes C\times W$$
where $C$ is the infinite cyclic group generated by $f_{\pi}$ and $W$ is
a direct product of countable many copies of
$\Z_p$.
Using these, we can easily find an open subgroup
$\N_w \subset \O_w^*$ so that $\O_w^*/\N_w \simeq
\Z/p^n\Z$, $f\in \N_w$ $\forall f \in \O_{S'}^*$
and $f_{\pi}$ generates a subgroup of order $p^{n-m_1}$ in
$\O_w^*/\N_w$.

\end{section}

\end{document}